\documentclass[11pt,leqno]{article}
\usepackage{enumerate}
\usepackage{amsmath}
\usepackage{amssymb}
\usepackage[english]{babel}
\usepackage{microtype}
\usepackage{authblk}

\usepackage{fullpage}

\newtheorem{thm}{Theorem}[section]
\newtheorem{lem}[thm]{Lemma}

\newtheorem{defi}[thm]{Definition}

\newcounter{claim}

\newenvironment{proof}[1][]%
 {\noindent {\setcounter{claim}{0}\sc proof ---
   }{#1}{}}{\hfill$\Box$\vspace{2ex}}

{\refstepcounter{claim}\vspace{1ex}\noindent{\it{\textit{Claim }\arabic{claim}{#1}{: }}}\it}{\vspace{1ex}}

	{\noindent {}{#1}{}}{ This proves Claim~\arabic{claim}.\vspace{1ex}}

\newcommand{\expe}{{\mathsf E}}
\newcommand{\prob}{{\mathsf{Pr}}}


\title{A graph inequality on the common neighbourhood}
\author[1]{Xiaomin Chen}
\author[2]{Fenglin Huang}
\author[2]{Shuhan Zhou}
\author[3]{\\ Mingxuan Zou}
\author[4]{Junchi Zuo}
\affil[1]{Shanghai Hypers Inc. (Shanghai, China)}
\affil[2]{Shanghai High School (Shanghai, China)}
\affil[3]{The High School Affiliated to Renmin University of China (Beijing, China)}
\affil[4]{No. 2 High School of East China Normal University (Shanghai, China)}

\begin{document}

\maketitle

\begin{abstract}
In this note we prove a graph inequality based on the sizes of the common neighbourhoods.
We also characterize the extremal graphs that achieve the equality.

The result was first discovered as a consequence of the classical Forster's theorem in
electric networks. 
We also present a short combinatorial proof that was inspired by a similar inequality related to 
the celebrated Tur\'an's theorem.
\end{abstract}

\section{Introduction}

A (simple) graph $G$ is a pair $G=(V, E)$ where $V$ is the vertex set and $E \subseteq \binom{V}{2}$ is the edge set. For any pair of distinct vertices $u$ and $v$, we denote $\{u, v\}$ by $uv$. For a vertex $u$, the (open) neighbourhood is defined as $N(u) = \{v \in V : uv \in E\}$. In particular, $u \not\in N(u)$. For two vertices $u$ and $v$, $N(u) \cap N(v)$ is the set of their common neighbours. For other standard graph notations not specified here, we refer to the textbook \cite{BM}.

In this short note we prove the following theorem for any connected graph.

\begin{thm}\label{thm.main}
Let $G=(V, E)$ be a connected graph with $n$ vertices, then 

(a).
\[
\sum_{uv \in E} \frac{1}{|N(u) \cap N(v)| + 2} \geq \frac{n-1}{2};
\]
(b). the equality in (a) holds if and only if every biconnected component of $G$ forms a clique.
\end{thm}

The theorem was first discovered and proved by one of the authors, Mingxuan, 
using electric networks. 
In Section \ref{sect.circuit} we give the original reasoning and a generalization.
While the implication from the classical results in electric networks is easy,
we find the theorem in its graph theoretical form in terms of common degrees is interesting.
In Section \ref{sect.graph} we give a simple combinatorial proof, 
which was inspired by one of the classical proofs of the classical Tur\'an's theorem.
In that proof we also give a cute lemma about the ordering of the vertices and common neighbours.

\section{Graph theoretical proofs}\label{sect.graph}

In one of the well known proofs (see \cite{AS}) of Tur\'an's theorem \cite{Turan}, we have the inequality 
(\cite{Caro} and \cite{Wei})
\[
\sum_{v \in V(G)} \frac{1}{|N(v)| + 1} \leq \alpha(G),
\]
where $\alpha(G)$ is the independence number of the graph $G$.
It was the similarity in the forms inspired us to prove Theorem \ref{thm.main} in a similar way.

\begin{defi} For a graph $G=(V, E)$, and any ordering of the vertices,
that is, an injection
$\pi : V \to \mathbb{N}$, call a pair $(u, v) \in V^2$ {\em good with respect to $\pi$} if $uv \in E$, $\pi(u) < \pi(v)$, and 
$\pi(u) < \pi(w)$ for any $w \in N(u)\cap N(v)$.

Define the graph $G^{(\pi)} = (V, E^{(\pi)})$ to be the subgraph whose edge set is formed by good pairs, that is,
\[
E^{(\pi)} = \{ \{u, v\} : (u, v) \mbox{ is good with respect to $\pi$}\}.
\]
\end{defi}

\begin{proof}[(of Theorem \ref{thm.main})]
Randomly uniformly pick an ordering $\pi : V \to [n]$, and define the random variables
$\chi_{u,v}$ to be the indicator of the event that $(u, v)$ is good, and $X = \sum_{(u,v), uv \in E} \chi_{u,v}$
be the number of the good pairs.
It is clear that, for any edge $uv$, 
\[
\expe \chi_{u,v} = \prob(\chi_{u,v} = 1) = \frac{1}{|N(u) \cap N(v)| + 2},
\]
so
\[
\expe X = \sum_{(u,v) \in V^2, uv \in E} \frac{1}{|N(u) \cap N(v)| + 2}.
\]
Note that each edge $uv$ contributes two terms in the sum. Thus, Theorem \ref{thm.main} follows directly from the next lemma.
\end{proof}

\begin{lem}\label{lem.main}
Let $G$ be a graph with $n$ vertices.

(a) $G^{(\pi)}$ is connected for any injection $\pi : V \to \mathbb{N}$.

(b) $G^{(\pi)}$ is a tree for every injection $\pi : V \to \mathbb{N}$ if and only if every biconnected component of $G$ forms a clique.
\end{lem}

We provide three short proofs to part (a) of Lemma \ref{lem.main} then a proof of part (b).

The first proof uses induction.

\begin{proof}[(Proof 1 to part (a))] We use an induction on $n$. The base case $n=1$ is obvious. 

Now suppose $G$ is a connected graph with $n>1$ vertices, $\pi$ is an ordering of the vertices, and $v$ is the last vertex in the ordering, that is, $\pi(v) = \max \{ \pi(u) : u \in V\}$.
Let $G_1, G_2, \dots, G_k$ ($k \geq 1$) be the connected components of the graph $G - v$, denote by $V_i$ the vertex set of $G_i$, and denote by $\pi_i$ the restriction of $\pi$ on $V_i$.

By the inductive hypothesis, $G_i^{(\pi_i)}$ is connected for every $1 \leq i \leq k$. 
It is easy to check that $G_i^{(\pi_i)}$ is a subgraph of $G^{(\pi)}$ ---
let $(a, b)$ be a good pair in $G_i$ with respect to $\pi_i$, we have $ab \in E(G_i) \subseteq E(G)$, $\pi(a) = \pi_i(a) < \pi_i(b) = \pi(b)$, and any vertex $c$ in $N_G(a) \cap N_G(b)$ is either $v$ or in the same connected component with $a, b$ in $G-v$, i.e. $c \in V_i$. In the former case $\pi(a) < \pi(c) = \pi(n)$; 
in the latter case $\pi(a) < \pi(c)$ since $c$ is also a common neighbour of $a$ and $b$ in $G_i$, and $(a, b)$ is good in $G_i$ with respect to $\pi_i$.

Let $a_i$ be the vertex in $N_G(v) \cap V_i$ with the smallest $\pi$-value. We claim that $(a_i, v)$ is 
also a good pair. Indeed, $a_i \in N_G(v)$ implies $a_i v \in E$; $\pi(a_i)<\pi(v)$ because $v$ has the maximum $\pi$-value, and any common neighbour $b$ of $a_i$ and $v$ is in $V_i$ and in $N_G(V)$,
so by the definition of $a_i$, we have $\pi(a_i) < \pi(b)$.

Thus, $G^{(\pi)}$ contains all the edges of $G_i^{(\pi_i)}$ that connected all the vertices in $V_i$,
and edges $va_i, va_j$ that connects $V_i$ and $V_j$, and $v$ and $V_i$'s; it is a connected graph. 
\end{proof}

The second proof shows that the edges of a ``minimum coded tree'' provides good pairs.

\begin{proof}[(Proof 2 to part (a))]
For any spanning tree $T$ of $G$, define the weight
\[
w(T) = \sum_{uv \in E(T)} \min(\pi(u), \pi(v)).
\]
Suppose $T^*$ is a tree with the minimum weight. 
We claim that any pair $(u, v)$ with $\pi(u) < \pi(v)$ and $uv \in E(T)$ is good with respect to $\pi$. 

Assume the contrary, there is $a \in N(u) \cap N(v)$ such that $\pi(a) < \pi(u)$. $T^* - uv$ has two connected components $C_u$ and $C_v$, where $u$ is in $C_u$ and $v$ is in $C_v$.
When $a$ is in $C_u$ (resp. $C_v$), $T' = T^* - uv + av$ (resp. $T^* - uv + au$) is a spanning tree of $G$,
but, in both cases,
\[
w(T') = w(T) - \pi(u) + \pi(a) < w(T),
\]
which is a contradiction.
\end{proof}

In the last proof we provide an algorithm that finds all the good pairs.

\begin{proof}[(Proof 3 to part (a))]
We maintain a graph $H$ and set $H = G$ in the beginning. For each edge $uv \in E$,
define $\sigma(uv) = \min(\pi(u), \pi(v))$. We do the following for each edge, in the order from
bigger $\sigma$ value to the smaller, and break ties arbitrarily: For each edge $uv$, if there
is a triangle $uva$ in $H$ with $\pi(a) < \pi(u)$ and $\pi(a) < \pi(v)$, we delete $uv$ from $H$.

On one hand, $H$ is always a subgraph of $G$. For any pair $(u, v)$ where $uv$ gets deleted,
the triangle $auv$ certifies that $(u, v)$ is not a good pair.
On the other hand, suppose $(u, v)$ is not a good pair, there is a triangle $auv$ in $G$ with $\pi(a) < \pi(u)$ and $\pi(a) < \pi(v)$. In the algorithm when we consider $uv$, by our ordering, none of $au$ and
$av$ is deleted from $H$, so $uv$ will be deleted.

Hence, the edge set of the final graph $H$ is indeed $G^{(\pi)}$. Note that we only delete an edge from a triangle from $H$, so $H$ is always connected.
\end{proof}

Now we prove the second part of Lemma \ref{lem.main}.

\begin{proof}[(part (b))] {\em If:} Let $G$ be a graph whose every biconnected component is a clique,
$\pi$ and ordering of its vertices. Assume the contrary that $G^{(\pi)}$ contains a cycle $C$ whose
vertices, sorted according to $\pi$, are
\[ \pi(v_1) < \pi(v_2) < \dots < \pi(v_k).\]
$v_2$ has two edges on $C$, so there is $i>2$ such that $(v_2, v_i)$ is a good pair. 
However, the cycle $C$ certifies that $v_1, v_2, v_i$ are in the same biconnected component,
which is a clique, so $v_1 \in N(v_2) \cap N(v_i)$ thus $(v_2, v_i)$ is not a good pair. A contradiction.

{\em Only if:}  Suppose $G$ has a biconnected component $H$ that is not a clique, 
so there are two points in $H$ with distance 2, therefore there are $a, b, c \in V(H)$ such that $ac, bc \in E$, $ab \not\in E$. Since $H$ is a biconnected component, there are paths between $a$ and $b$
that do not go through $c$. Let $P = (a, v_1, v_2, \dots, v_k, b)$ be such a path with smallest $k$.
We have $k \geq 1$.

Now let $\pi$ be a permutation of the vertices where
\[
\pi(a) < \pi(b) < \pi(c) < \pi(v_1) < \pi(v_2) < \dots < \pi(v_k),
\]
and $\pi(v_k) < \pi(u)$ for any $u \not\in \{a, b, c, v_1, \dots, v_k\}$.

It is clear that $(a, c)$ and $(a, v_k)$ are good pairs. $(b, c)$ and $(b, v_1)$ are also good pairs because
the only vertex appears before $b$ is $a$ and it is not connected to $b$. Next we prove that $(v_i, v_{i+1})$ is a good pair for any $1 \leq i < k$ --- otherwise, there is $x$ such that $\pi(x) < \pi(v_i)$
and $xv_iv_{i+1}$ is a triangle; but, by the minimality of $P$, 
$x \neq a$ since $av_i \not\in E$; 
$x \neq b$ since $bv_{i+1} \not\in E$;
$x \neq v_j$ for any $j<i$ since $v_jv_{i+1} \not\in E$.

So, all the edges on $P$ are in $G^{(\pi)}$, and $ab$, $ac$ are also in $G^{(\pi)}$; they form a cycle.
\end{proof}
 
\section{Proof via the electric networks}\label{sect.circuit}

Given a connected graph $G=(V, E)$, we construct an electric network $\mathcal{N}_G$ on the set of vertices,
and connect each pair of adjacent vertices by a resistance of magnitude 1. The resistance distance
$R_{u,v}$ between two vertices is the effective resistance between $u$ and $v$ in $\mathcal{N}_G$. It is well known that this is a metric on $V$.

Theorem \ref{thm.main} was first discovered by Mingxuan while pondering over
the following beautiful fact in electric networks, known as Forster's first theorem. 
(See \cite{Forster}, \cite{Klein}.)

\begin{thm}\label{lem.resist} \cite{Forster}
For a connected graph $G=(V, E)$, denote $R_{u,v}$ the effective resistance between two vertices $u$
and $v$ in $\mathcal{N}_G$, then the sum of effective resistances over the edges
\[
\sum_{uv \in E} R_{u,v} = n-1.
\]
\end{thm}

Theorem \ref{thm.main} is a consequence of Forster's theorem.

\begin{proof}[(of Theorem \ref{thm.main})]

(a). 
For any edge $uv \in E$, consider the network modified from $\mathcal{N}_G$ by keeping the resistances on edges
\[
\{ uv \} \cup \{uw, vw : w \in N(u) \cap N(v)\},
\]
and deleting all the other edges (or, equivalently, we increase all the other resistances towards $+\infty$), 
the resulting network is equivalent to $|N(u) \cap N(v)|+1$ parallel edges, one having resistance 1 and others
having resistance 2 each.
By Rayleigh’s monotonicity law,
the effective resistance $R_{u,v}$ is upper-bounded by
\begin{equation}\label{eq.resist_bound}
R_{u,v} \leq \frac{1}{\frac{1}{2} |N(u) \cap N(v)| + 1}.
\end{equation}
So, combined with Lemma \ref{lem.resist},
\[
\sum_{uv \in E} \frac{1}{|N(u) \cap N(v)| + 2} = \frac{1}{2}\sum_{uv \in E} \frac{1}{\frac{1}{2} |N(u) \cap N(v)| + 1} 
\geq \frac{1}{2} \sum_{uv \in E} R_{u,v} = \frac{n-1}{2}.
\]

(b). It is easy to check the {\em if} direction. To check the {\em only if} direction, we assume some biconnected component $C$ of $G$ is not a clique, and prove that for at least one pair of adjacent vertices the equality in \eqref{eq.resist_bound} does not hold.

Being a connected, non-complete graph, a well-known fact is that $C$ has an induced $P_3$,
i.e., $u, v, z \in V(C)$ such that $uv, uz \in E$, $vz \not\in E$. Since $C$ is biconnected, $C-u$
is connected, pick a shortest path $P$ from $z$ to $\{v\} \cup (N(u) \cap N(v))$ in $C-u$.
Now, in $\mathcal{N}_G$, only keep the edges
\[
\{ uv \} \cup \{uw, vw : w \in N(u) \cap N(v)\} \cup P.
\]
It is easy to see 
\[
R_{u,v} < \frac{1}{\frac{1}{2} |N(u) \cap N(v)| + 1}.
\]
\end{proof}

From here it is easy to get a generalization of Theorem \ref{thm.main}. For any pair of vertices $u$ and $v$, define $P(uv, \ell)$ to be the maximum number of internally vertex-disjoint paths between $u$ and $v$ 
with length between 2 and $\ell$, inclusively, then the effective resistance
\[
R_{u,v} \leq \frac{1}{1+P(uv,\ell)/\ell}.
\]
Thus,

\begin{thm}
In a connected graph $G=(V, E)$ with $n$ vertices, 
\[
\sum_{uv \in E} \frac{1}{P(uv,\ell) + \ell} \geq \frac{n-1}{\ell}.
\]
\end{thm}

\section*{Acknowledgement}

Theorem \ref{thm.main} was discovered by Mingxuan and proved by him using the method in Section \ref{sect.circuit}. It was provided to the Semicircle math club as one of the challenges whether it has a simple combinatorial proof, and we discovered the proofs in Section \ref{sect.graph} during our club discussions.

We thank the members of the Semicircle Club for the pleasant discussions. And we thank Mario Szegedy
for helpful discussions and comments.

\end{document}